\documentstyle[11pt]{article}

\input{mssymb}
\def\german{\frak}

\newtheorem{thm}{Theorem}[section]

\newtheorem{theorem}[thm]{Theorem}

\newtheorem{proposition}[thm]{Proposition}

\newtheorem{cor}[thm]{Corollary}

\newtheorem{lemma}[thm]{Lemma}

\newenvironment{example}{\begin{example2} \rm}{\end{example2}}
\newtheorem{example2}[thm]{Example}

\newenvironment{remark}{\begin{remark2} \rm}{\end{remark2}}
\newtheorem{remark2}[thm]{Remark}

\newenvironment{definition}{\begin{defn2} \rm}{\end{defn2}}
\newtheorem{defn2}[thm]{Definition}

\newcommand{\pr} {\smallskip\noindent{\sl Proof.\,\,}}
\newcommand{\spr} {\smallskip\noindent{\sl Sketch of proof.\,\,}}

\newenvironment{proof}	{\pr}{\hspace*{\fill}$\Box$\\}

\newcommand{\bbR}{{\Bbb R}}

\newcommand{\bbC}{{\Bbb C}}

\def\Phinv{\Phi^{-1}}
\newcommand\inv {{^{-1}}}
\newcommand\calV {{\cal V}}

\def\g{{\german g}}  \def\fg{{\german g}} 
\def\t{{\german t}}  \def\ft{{\german t}}
  
\def\z{{\german z}}   
\def\h{{\german h}}   \def\fh{{\german h}}
\def\fm{{\german m}}
 
\def\cpwc{\ft_+^*}

\def\tU{{\tilde{U}}}

\def\ra{\rightarrow}
\def\ol{\overline}

\begin{document}

\title{Non-abelian convexity by symplectic cuts\\ dg-ga/9603015}

\author{Eugene Lerman\thanks{Partially supported by a National Science
Foundation Postdoctoral Fellowship.}\\ {\small Department of
Mathematics, Univ.\ of Illinois, Urbana, IL 61801\
lerman@math.uiuc.edu}\\[2pt] Eckhard Meinrenken,\thanks{Supported by a
Feodor-Lynen Fellowship of the Humboldt Foundation.} \, Sue
Tolman,\thanks{Partially supported by a National Science Foundation
Postdoctoral Fellowship.} \,and Chris Woodward\thanks{Supported by a
Sloan Foundation Doctoral Dissertation Fellowship.}  \\ \small
Department of Mathematics, M.I.T., Cambridge MA 02139, \\ \small
mein@math.mit.edu, tolman@math.mit.edu, woodward@math.mit.edu}

\maketitle

\begin{abstract}
In this paper we extend the results of Kirwan et alii on convexity 
properties of the moment map for Hamiltonian group actions,
and on the connectedness of the fibers of the moment map, 
to the case of non-compact orbifolds. 

Our motivation is twofold.  First, the category of orbifolds is
important in symplectic geometry because, generically, the symplectic
quotient of a symplectic manifold is an orbifold.  Second, our proof
is conceptually very simple since it reduces the non-abelian case to
the abelian case.
\end{abstract}

\tableofcontents

\section{Introduction}

In this paper we prove the following theorem, which extends a result
of Kirwan (\cite{Kconvex}, \cite{Kbook}) to the case of orbifolds
which need not be compact.  Recall that a subset $\Delta$ of a vector
space $V$ is {\bf polyhedral} if it is the intersection of finitely
many closed half-spaces.  Recall also that $\Delta \subset V$ is {\bf
locally polyhedral} if for any point $x\in \Delta$ there is a
neighborhood $U$ of $x$ in $V$ and a polyhedral set $P$ in $V$ such
that $U\cap \Delta = U\cap P$.

\begin{thm} \label{MainThm} 
Let $(M,\omega)$ be a connected symplectic orbifold with a Hamiltonian
action of a compact Lie group $G$ and a proper moment map $\Phi:M \to \g^*$.
\begin{enumerate}
\item \label{ThmA} Let $\cpwc$ be a closed Weyl chamber for the Lie
group $G$ considered as a subset of $\g^*$.  The {\em moment set}
$\Phi(M) \cap \cpwc$ is a convex locally polyhedral set. In
particular, if $M$ is compact then the moment set is a convex
polytope.
\item \label{ThmB}
Each fiber of the moment map $\Phi$ is connected.
\end{enumerate}
\end{thm}

Kirwan's theorem has numerous applications in symplectic geometry:  For
example, it is used in the classification of Hamiltonian $G$-orbifolds, in
geometric quantization, and in the study of the existence of invariant
K\"ahler structures.

Convexity theorems in symplectic geometry have a long history.  For
the action of a maximal torus on a coadjoint orbit Theorem
\ref{MainThm} was proved by Kostant \cite{Kostant}, extending previous
results of Schur and Horn.  This was generalized to the actions of
subgroups by Heckman \cite{Hec}. For Hamiltonian torus actions on
manifolds the theorem was proved independently by Atiyah \cite{A} and
by Guillemin and Sternberg \cite{GSc}. In the non-abelian case
Guillemin and Sternberg proved that the moment set is a union of
convex polytopes and that it is a single polytope for a K\"ahler
manifold. The convexity in the K\"ahler case was independently proved by
Mumford \cite{NM}. The first complete proof for Hamiltonian actions of
non-abelian groups on manifolds was given by Kirwan \cite{Kconvex},
using Morse theory of the Yang-Mills functional and results of
Guillemin and Sternberg.  Several alternative proofs of Kirwan's
result have appeared: \cite{CDM}, \cite{Hilgert} (based on \cite{CDM})
and \cite{Sj} (based on \cite{B}).

Our motivation for offering this proof is twofold.  First,
Theorem~\ref{MainThm} generalizes the result from manifolds to
orbifolds; this category is important in symplectic geometry because,
generically, the symplectic quotient of a symplectic manifold is an
orbifold.  Secondly, our proof is conceptually very simple.  The first
step is to generalize the abelian version of the theorem to certain
non-compact orbifolds. The key idea is to apply the technique of
symplectic cutting \cite{L} to reduce to the compact case.  In the
second step we reduce the non-abelian case to the abelian case by
means of a symplectic cross-section: We show that there is a unique
open wall $\sigma$ of the Weyl chamber $\cpwc$ such that
$\Phi(M)\cap\sigma$ is dense in the moment set
$\Delta=\Phi(M)\cap\cpwc$, and that the preimage $\Phi^{-1}(\sigma)$
is a $T$-invariant symplectic suborbifold.

Under assumptions similar to those of the main theorem, the result
holds for Hamiltonian actions of loop groups provided the preimage of
an alcove under the moment map is finite dimensional. This will be
discussed in a forthcoming paper.

In the last section, using similar methods we extend to symplectic
orbifolds a result of Sjamaar \cite{Sj} which says that the the moment
set near a point $x$ can be read off from local data near an arbitrary
point in the fiber $\Phinv(x)$:

\begin{thm}
Let $(M,\omega)$ be a connected symplectic orbifold with a Hamiltonian
action of a compact Lie group $G$ and a proper moment map $\Phi:M \to \g^*$.
For every  $m\in \Phi^{-1}(\cpwc)$, and 
every $G$-invariant neighborhood  $U$ of $m$ in $M$, the image 
$\Phi(U)\cap \cpwc$ is a relatively open neighborhood of $\Phi(m)$ in 
$\Phi(M)\cap\cpwc$.
\end{thm}

\section{Background}

In this section, we extend some well-known results in symplectic
geometry from manifolds to or\-bi\-folds.  These extensions are
straightforward.  Readers who are familiar with the standard versions
of these results and not too skeptical may wish to skip this section.

\subsection{Hamiltonian actions on orbifolds}
\label{secab}

An {\bf orbifold} $M$ is a topological space $|M|$, together with an
{\bf atlas} of {\bf uniformizing charts} $(\tU,\Gamma,\varphi)$, where
$\tU$ is open subset of $\bbR^n$, $\varphi(\tU)$ is an open subset
of $|M|$, $\Gamma$ is a finite group which acts linearly on $\tU$ and
fixes a set of codimension at least two, and $\varphi: \tU \to |M|$
induces a homeomorphism from $\tU/\Gamma$ to $\varphi(\tU) \subset
|M|$.  Just as for manifolds, these charts must cover $|M|$; they are
subject to certain compatibility conditions; and there is a notion of
when two atlases of charts are equivalent.  For more details, see
Satake \cite{Sa}.

A smooth function on $M$ is a collection of smooth invariant functions
on each uniformizing chart $(\tU,\Gamma,\varphi)$ which agree on
overlaps of the images $\varphi (\tU)$. Differential forms, vectors
fields, and other objects can be similarly defined.  There is also a
notion of morphisms (maps) of orbifolds.

Let $x$ be a point in an orbifold $M$, and let $(\tU,\Gamma, \varphi)$
be a uniformizing chart with $x\in\tU/\Gamma$.  The {\bf
(orbifold) structure} group of $x$ is the isotropy group $\Gamma_x$ of
$\tilde{x} \in \tU$, where $\varphi(\tilde{x}) = x$.  The group
$\Gamma_x$ is well defined as an abstract group. The  
tangent space to $\tilde{x}$ in $\tU$, considered as a representation 
of  $\Gamma_x$ is called the {\bf uniformized tangent space at $x$}, 
and denoted  by $\tilde{T}_xM$. 
The quotient $\tilde{T}_x M/\Gamma_x$ is
$T_x M$, the fiber of the {\bf tangent bundle} of $M$ at $x$.

Let $G$ be a compact connected Lie group with Lie algebra $\g$.  A
{\bf smooth action} $a$ of $G$ on an orbifold $M$ is a smooth orbifold
map $a: G\times M \to M$ satisfying the usual laws for an action.
Given an action of $G$ on $M$, every vector $\xi \in \g$ induces a
vector field $\xi_M$ on $M$.

A {\bf symplectic orbifold} is an orbifold $M$ with a closed
non-degenerate $2$-form $\omega$.  A group $G$ acts {\bf
symplectically} on $(M,\omega)$ if the action preserves $\omega$.  A
{\bf moment map } for a symplectic action of a group $G$ is an
equivariant map $\Phi: M \to \g^*$ such that
\begin{equation}\label{moment map defining equation}
i_{\xi_M} \omega = -\,\langle\, d\Phi,\xi
\rangle\quad \mbox{ for all }\quad\xi \in \g.
\end{equation}
If a moment map exists, we say that the action of $G$ on $(M,\omega)$
is Hamiltonian and refer to $(M,\omega)$ as a {\bf Hamiltonian
$G$-space}.  If $b \in \g^*$ is a regular value of the moment map
$\Phi$, then the orbifold version of the Marsden-Weinstein-Meyer
theorem says that the quotient $M_b := \Phinv(b)/G_b$ is a symplectic
orbifold, called the {\bf symplectic reduction} of $M$ at $b$.  
For a proof, see \cite{LT}.

The following theorem is due in the manifold case to Atiyah \cite{A}
and to Guillemin and Sternberg \cite{GSc}.  The proofs in the orbifold
case are in \cite{LT}.  Let $T$ be a torus, $\ft$ its Lie algebra,
$\ft^*$ the dual of $\ft$ and $\ell = \ker\{ \exp : \ft \to T\}$ the
{\bf integral lattice}.  A polytope in $\ft^*$ is {\bf rational} if
the faces of the polytope are cut out by hyperplanes whose normal
vectors are in the lattice $\ell$.
\begin{theorem} \label{AbConv} 
Let $(M,\omega)$ be a compact, connected symplectic orbifold, and let
$\Phi: M \to \t^*$ be a moment map for a Hamiltonian torus action on
$M$.
\begin{enumerate}
\item The image $\Phi(M)$ is a rational convex polytope, and
\item each fiber of $\Phi$ is connected.
\end{enumerate}
\end{theorem}

\subsection{Symplectic Cuts} \label{Symplectic Cutting}

Symplectic cutting is a technique which allows one to naturally
construct a symplectic structure on a subquotient of a symplectic
orbifold.  It was introduced in \cite{L}.  We will use symplectic cuts
to produce compact orbifolds out of non-compact ones.

Let $(M,\omega)$ be a symplectic orbifold with a Hamiltonian circle action
and a moment map $\mu :M \to \bbR$.  Suppose that $\epsilon$ is a
regular value of the moment map. 
Consider the disjoint union 
$$ M_{[\epsilon,\infty)}:=\mu^{-1}((\epsilon,\infty))\cup
M_\epsilon,$$ 
obtained from the orbifold with boundary $\mu^{-1}([\epsilon,\infty))$
by collapsing the boundary under the $S^1$-action. We claim 
that $ M_{[\epsilon,\infty)}$  
admits a natural structure of a symplectic orbifold, in such a way that 
the embeddings of $\mu^{-1}((\epsilon,\infty))$ and $M_\epsilon$ are 
symplectic.
Moreover, the induced circle action on $M_{[\epsilon,\infty)}$ 
is Hamiltonian, with
a moment map coming from the restriction of the original moment map
$\mu$ to $\{m\in M: \mu (m) \geq \epsilon\}$.

\begin{definition}  
We call the symplectic orbifold $M_{[\epsilon, \infty)}$ the {\bf
symplectic cut } of $M$ with respect to the ray $[\epsilon, \infty)$
(with symplectic form and moment map understood).
\end{definition}

To see that the claim holds consider the product $M \times \bbC$ of
the orbifold with a complex plane.  It has a natural (product)
symplectic structure : $\omega + (-i) dz \wedge d\bar{z}$.  The
function $\nu: M \times \bbC \to \bbR$ given by $ \nu(m,z) = \mu(m) -
| z |^2$ is a moment map for the diagonal action of the circle (it
commutes with the original action of the circle on $M$).  The point
$\epsilon$ is a regular value of $\mu$ if and only if it is a regular
value of $\nu$.  The map
$$ 
\{m\in M :\mu \geq \epsilon\} \to \nu ^{-1} (\epsilon), 
\qquad m\mapsto (m,\sqrt{ \mu (m) - \epsilon}) 
$$ 
descends to a homeomorphism from the cut space $M_{[\epsilon, \infty)}$
to the reduced space $\nu ^{-1} (\epsilon)/S^1$. Moreover the 
homeomorphism is $S^1$-equivariant and it is a symplectomorphism on
the set $\{\mu > \epsilon\}$.

\begin{example}  Consider the complex plane $\bbC$ with its standard 
symplectic form $(-i) dz \wedge d\bar{z}$, and let a circle $U(1) =
\{z\in \bbC: |z| =1 \}$ act by multiplication.  The moment map is $\mu
(z)= -|z|^2$.  The symplectic cut of $\bbC$ at $\epsilon <0$ is a
two-sphere. 

\end{example}

This construction generalizes to torus actions as
follows.\footnote{See \cite{W}, \cite{Mein} for the generalization to
non-abelian actions.}  Let $\mu:M \to \t^*$ be a
moment map for an { effective }action of a torus $T$ on a
symplectic orbifold $(M,\omega)$ and let $\ell \subset \t$ denote the
integral lattice.  
Choose $N$ vectors $v_j \in \ell$, $1 \leq j \leq
N$. The form $\omega -i \sum_{j=1}^N\, dz_j\wedge d \bar{z}_j$ is a
symplectic form on the orbifold $M \times \bbC^N$.  The map $\nu:M
\times \bbC^N \to \bbR^N$ with $j$th component 
$$ 
\nu_j(m,z) = \langle\mu(m), v_j \rangle - | z_j |^2 
$$ 
is a moment map for a $(S^1)^N$ action on $M \times \bbC^N$.  For any $b
\in \bbR^N$, define a convex rational polyhedral set
$$
P = \{ x\in \t^* \mid \langle x, v_i \rangle \geq b_i
\mbox{ for all } 1 \leq i \leq N \}.
$$ 
The {\bf symplectic cut of $M$ with respect to a rational polyhedral
set $P$} is the reduction of $M \times \bbC^N$ at $b$. We denote it by
$M_{P}$.

\begin{remark}  \label{Generic}  
If $b$ is a regular value of $\nu$, then $M_{P}$ is a symplectic
orbifold.  Note that regular values are generic.  Thus, for an
intersection $P$ of finitely many rational half-spaces of $\t^*$ in
general position the cut space $M_P$ is a symplectic orbifold. Note
that if $P$ is a compact polytope, then the fact that $P$ is generic
implies that $P$ is {\bf simple:} the number of codimension one faces
meeting at a given vertex is the same as the dimension of $P$.
\end{remark}  

There is a natural decomposition of
$M_{P}$ into a union of symplectic suborbifolds, indexed by the open
faces $F$ of $P$: 
$$ 
M_{P}=\bigcup_F \,\mu^{-1}(F)/T_F, 
$$ 
where $T_F\subset T$ is the torus {\bf perpendicular} to $F$.  That
is, the Lie algebra of $T_F$ is the annihilator of the linear space
defined by the face $F$.  Thus, one can think of $M_{P}$ as
$\mu^{-1}(P)$, with its boundary collapsed by means of the
$T_F$-actions.  Under this identification, the $T$-action on $M$
descends to a $T$-action on $M_P$. One can show that this action is
smooth and Hamiltonian, with moment map $\mu_P$ induced from the
restriction $\mu|\mu^{-1}(P)$.  In particular,
$$
\mu_P(M_{P})=\mu(M)\cap P.
$$ 
It follows immediately that the cut space $M_{P}$ is compact exactly if
$\mu^{-1}(P)$ is.  Similarly $M_{P}$ is connected if and only if
$\mu^{-1}(P)$ is connected.  To summarize, we have the following result.

\begin{proposition}\label{prop_symplectic_cuts}
Let $\mu:M \to \t^*$ be a moment map for an {\em effective }action of
a torus $T$ on a symplectic orbifold $(M,\omega)$.  Let $P \subset
\ft^*$ be a generic rational polyhedral set. Let ${\cal F}$ be the set
of all open faces of $P$.  Then the topological space $M_P$ defined by
$$
M_{P}=\bigcup_{F\in {\cal F}} \,\mu^{-1}(F)/T_F,
$$
where $T_F\subset T$ is the subtorus of $T$ perpendicular to $F$, is a
symplectic orbifold with a natural Hamiltonian action of the torus
$T$.  Moreover, the map $\mu _P: M_P \to \ft^*$, induced by the
restriction $\mu |\mu\inv (P)$, is a moment map for this action.
Consequently, 
\begin{enumerate}
\item the cut space $M_P$ is connected if and only if $\mu \inv (P)$
is connected;

\item the fibers of $\mu_P$ are connected if and only if the fibers of
$\mu| \mu \inv (P)$ are connected;

\item $M_P$ is compact if and only if $\mu \inv (P)$ is compact.

\end{enumerate}
\end{proposition}

\subsection{Principal orbit type}
Consider a connected orbifold $M$, together with an action of 
a compact connected Lie group $G$. For each $m\in M$, let 
$\g_m\subset\g$ be the corresponding isotropy Lie algebra. Clearly, 
$\g_{g\cdot m}={\rm Ad}(g)(\g_m)$ for all $g\in G$. The set of subalgebras
$$
(\g_m) = \{{\rm Ad}(g)(\g_m) \mid g\in G \}
$$ 
is called the (infinitesimal) {\bf orbit type} of $m$. As in the case
of manifolds, each point $m\in M$ has a neighborhood $U$ such that the
number of orbit types $(\g_{m'})$, for $m'\in U$ is finite and each
$\g_{m'}$ is subconjugate to $\g_m$.

\begin{proposition}
There exists a unique orbit type $(\h)$ (called the {\bf principal
orbit type}) such that the set
$$ M_{(\h)}=\{m\in M|\,(\g_m)=(\h)\} $$
of points of orbit type $(\h)$ is open, dense and connected.
\end{proposition}
\begin{proof}
The proof is an adaptation of the proof for manifold to orbifolds.
The key point is that slices (Definition~\ref{definition of slice})
exist for actions of compact groups on orbifolds \cite{LT}.
\end{proof}
 
\begin{definition}
For an action of a compact connected Lie group $G$ on a connected 
orbifold $M$, we define the {\bf principal stratum} 
$M_{prin}$ to be the intersection of the set 
$M_{(\h)}$ of points of principal orbit type with the set $M_{smooth}$
of smooth points  of $M$. 
\end{definition}

\begin{remark}
By definition, the set $M_{sing}$ of singular points of an orbifold
$M$ is a union of submanifolds of codimension at least 2. Therefore
$M_{smooth}$ is open, dense and connected, and so is $M_{prin}$.
\end{remark}

\begin{remark}\label{cons_rank_remark}
Let $M$ be a connected Hamiltonian $G$-orbifold. As in the case of
manifolds, the definition of the moment map, equation (\ref{moment map
defining equation}), implies that the image of the differential of the
moment map at a point $m$ is the annihilator of the corresponding
isotropy Lie algebra $\fg_m$. Consequently the restriction of the
moment map to $M_{prin}$ has constant rank.
\end{remark}

\section{The principal cross-section}

In this section, we define cross-sections and show that they have the
properties needed to reduce the non-abelian case to the abelian case.

\subsection{The cross-section theorem}

Theorem~\ref{cross-section thm} is a version of the cross section
theorem of Guillemin and Sternberg \cite[Theorem~26.7]{GS} adapted to
orbifolds. This version of the theorem and its relationship to
fibrations of coadjoint orbits and representation theory are discussed
in \cite{GLS}. Since \cite{GLS} is not yet published, we thought that
it would be useful to present the proof.  

\begin{definition}\label{definition of slice}
Suppose that a group $G$ acts on an orbifold $M$.  Given a point $m$
in $M$ with isotropy group $G_m$, a suborbifold $U\subset M$
containing $m$ is a {\bf slice} at $m$ if $U$ is $G_m$-invariant,
$G\cdot U$ is a neighborhood of $m$, and the map
$$ 
G\times _{G_m} U \to G\cdot U, \quad [a, u]\mapsto a\cdot u 
$$ 
is an isomorphism.  In other words, $G\cdot y \cap U = G_m
\cdot y$ and $G_y \subset G_m$ for all $y\in U$. 
\end{definition}

\begin{remark}\label{remark natural slice}
Consider the coadjoint action of a compact connected Lie group $G$ on
$\g^*$. For all $x\in\g^*$, there is a unique largest open subset 
$U_x\subset\g_x^*\subset\g^*$ which is a slice at $x$. We refer to 
$U_x$ as the {\bf natural slice } at $x$ for the coadjoint action.
In order to describe the natural slice, we may assume without 
loss of generality that $x\in\cpwc$.
Let $\tau\subset\cpwc$ be the open
wall of $\cpwc$ containing $x$ and let $G_\tau $ denote the isotropy
Lie group of $x$ (the group is the same for all points of $\tau$). Then
$$ 
	U_\tau = G_\tau \cdot \{ y \in \cpwc \mid G_y \subset G_\tau\}
               = G_\tau \cdot \bigcup_{\tau\subset\ol{\tau'}}\tau'  
$$
is an open subset of $\fg_\tau ^*$, and is equal to the natural 
slice $U_x$.
\end{remark}

\begin{theorem}[Cross-section]\label{cross-section thm}
Let $(M, \omega)$ be a symplectic orbifold with a moment map $\Phi: M
\to \fg^*$ arising from an action of a compact Lie group $G$.  Let $x$
be a point in $\fg^*$ and let $U$ be the natural slice at $x$ (see
above).  Then the {\bf cross-section} $R:=\Phi ^{-1} (U)$ is a
$G_x$-invariant symplectic suborbifold of $M$, where $G_x$ is the
isotropy group of $x$.  Furthermore the restriction $\Phi |_R$ is a
moment map for the action of $G_x$ on $R$.
\end{theorem}

\begin{proof}
By definition of the slice, coadjoint orbits intersect $U$
transversally.  Since the moment map is equivariant, it is transversal
to $U$ as well.  Hence the cross-section $R=\Phi ^{-1} (U)$ is a
sub\-or\-bi\-fold.  Since the slice $U$ is preserved by the action of $G_x$
and the moment map is equivariant, the cross-section is preserved
by $G_x$.
It remains to show that for all $r\in R$, the uniformized tangent space
$\tilde{T}_r R$ is a symplectic subspace of $\tilde{T}_r M$.
Let $y=\Phi(r)$, and 
let $\fm$ be the $G_x$-invariant complement of $\fg_x$ in $\fg$.  Then
$T_yU$ is the annihilator of $\fm$ for any $y\in U$.  We claim that
\begin{itemize}
\item[(a)]  the tangent space $\tilde{T}_r R$ is symplectically 
perpendicular to $\fm_M (r) = \{\xi _M (r): \xi \in \fm\}$, the
subspace of the tangent space to the orbit through $r$ spanned by
$\fm$;
\item[(b)] $\fm_M (r)$ is a symplectic subspace of  $\tilde{T}_r M$.
\end{itemize}
Together the two claims establish the theorem because $\tilde{T}_r M =
\tilde{T}_r R \oplus \fm_M (r)$.  To see that (a) holds observe that
for $v\in \tilde{T}_r R$ and $\xi \in \fm$,
$$
\omega_r (\xi _M(r), v) = \langle \xi , d\Phi _r (v)\rangle = 0
$$
since $d\Phi _r (v) \in T_y U = \fm^\circ$ and $\xi \in \fm$.

To see that (b) is true,  observe that for $\xi, \eta \in \fg$,
$$
\omega_r (\xi _M(r), \eta _M (r)) = \langle \xi , d\Phi _r (\eta _M
(r))\rangle =  \langle \xi, ad^\dagger (\eta)\cdot \Phi (r)\rangle =
\langle [\xi, \eta], y\rangle.
$$
Thus $\fm_M (r)$ is symplectic if and only if $ad^\dagger (\fm) y$ is
a symplectic subspace of the tangent space $T_y (G\cdot y)$.
Since $G_x \cdot y \subset U$ and since $\fm = (T_y U)^\circ$, 
for any $\xi \in \fm$ and any $\zeta \in \fg_x$ we have
$$
\langle [\xi, \zeta], y\rangle =\langle \xi, ad^\dagger (\zeta)
y\rangle = 0,  
$$ 
i.e., $T_y (G_x \cdot y)$ and $ad^\dagger (\fm) y$ are symplectically
perpendicular in $T_y (G\cdot y)$.  Since $T_y (G\cdot y)=T_y (G_x
\cdot y)\oplus ad^\dagger (\fm) y $, it remains to show that the orbit
$G_x \cdot y$ is a symplectic submanifold of the coadjoint orbit
$G\cdot y$.  Since the natural projection $\pi : \fg^* \to \fg_x^*$ is
$G_x$-equivariant, $\pi (G_x \cdot y) = G_x \cdot \pi (y)$.  By the
definition of the symplectic forms on a coadjoint orbit the
restriction of the symplectic form on $G\cdot y$ to $G_x \cdot y$ is
the pull-back by $\pi$ of the symplectic form on the $G_x$ coadjoint
orbit $G_x \cdot \pi (y)$.
\end{proof}

\begin{remark}\label{assumptions}
For the theorem to hold for a non-compact Lie group $G$ one has to
assume that a slice $U$ exists at $x$ and that the differential of the
restriction to $U$ of the projection $\fg^* \to \fg_x^*$ is
surjective.
\end{remark}

\begin{remark}\label{bundle_over_coadj_orbit}
The cross-section $R$ need not be a slice for an action of $G$ on $M$
since the group $G_x$ need not appear as an isotropy group of any
point in the cross-section.  However, the set $G\cdot R$ of orbits
through the cross-section is an open subset of the manifold $M$ and it
is equivariantly diffeomorphic to the associated bundle
$G\times_{G_x}R$ over the coadjoint orbit $G\cdot x$.  This is because
$U$ is a slice at $x\in \fg^*$, $R= \Phi\inv (U)$ and $\Phi$ is
equivariant.
\end{remark}

\begin{remark}\label{general argument} 
Let $G$ be a compact connected group, and $M$ a Hamiltonian
$G$-orbifold, with moment map $\Phi:\,M\rightarrow \g^*$. 
In various applications in this paper, we will use
symplectic cross-sections to reduce statements about
$G$-orbits in $M$ to the case that the orbit is contained in the zero
level set $\Phinv(0)$.  The general argument is as follows. Let
$m\in\Phinv(\cpwc)$, and let $\tau\subset\cpwc$ be the open wall
containing $x=\Phi(m)$.  Let $U_\tau\subset\g^*$ be the natural slice and
$R_\tau$ the corresponding natural cross-section, which is a
Hamiltonian $G_\tau$-space.  Since $G_\tau$ contains the maximal
torus, there is a unique $G_\tau$-invariant decomposition
$$ \g=\z(\g_\tau)\oplus [\g_\tau,\g_\tau]\oplus{\frak m}_\tau,$$
where  $\z(\g_\tau)$ is the center of $\g_\tau$, 
$[\g_\tau,\g_\tau]$ its semi-simple
part and ${\frak m}_\tau$ a complement in $\g$.
Notice that $\z(\g_\tau)$ can be characterized as the fixed point 
set of the $G_\tau$-action on $\g$. It follows that the linear 
span of $\tau$ is equal to $\z(\g_\tau)^*$. 
Since $x=(\Phi|R_\tau)(m)\in \tau\subset\z(\g_\tau)^* $, one can shift the 
moment map $\Phi|R_\tau$ by $x$ to obtain a new moment map 
$\Phi_\tau'$ for the $G_\tau$-action on $R_\tau$
for which $m\in (\Phi_\tau')^{-1}(0)$.  
\end{remark}

\subsection{The principal wall and the corresponding cross-section}
The main result of this subsection is: 
\begin{theorem}\label{princ cross-section}
Let $G$ be a compact connected Lie group, and 
$M$ a connected Hamiltonian $G$-orbifold, with moment
map $\Phi: M\rightarrow\g^*$.
\begin{enumerate}
\item
There exists a unique open wall $\sigma$ of the  Weyl chamber 
$\cpwc$ with 
the property that $\Phi(M)\cap\sigma$ is dense in 
$\Phi(M)\cap\cpwc$. 
\item The preimage $Y=\Phi\inv(\sigma)$ is a connected 
symplectic $T$-invariant suborbifold of $M$, 
and the restriction $\Phi_Y$ of $\Phi$ to $Y$ is a moment
map for action of the maximal torus $T$.  
\item The set $G\cdot Y = \{g\cdot m \mid g\in G, \,\, m\in Y\}$ is
dense in $M$. 
\end{enumerate}
\end{theorem}

\noindent
We refer to $\sigma$ as the {\bf principal wall} and to
$Y=\Phi^{-1}(\sigma)$ as the {\bf principal cross-section}.
Lemma~\ref{lemma stabilizer} below is used to show the existence of
the principal wall.

\begin{lemma}\label{lemma stabilizer}
Let $G$ be a compact, connected Lie group, and $M$ a connected
Hamiltonian $G$-orbifold, with moment map $\Phi:M\rightarrow\g^*$. 
\begin{enumerate}
\item For all $m\in M_{prin}$, the isotropy Lie algebra $\g_m$ is an
ideal in $\g_{\Phi(m)}$, i.e., $[\fg_m, \fg_{\Phi(m)}] \subset
\fg_{m}$.
\item 
All points in the intersection $\Phinv(\cpwc)\cap M_{prin}$ 
have the same isotropy Lie algebra 
$\h\subset\g$. 
\item 
Given $\alpha\in\Phi(M_{prin})\cap\cpwc$ let $S$ be the affine
subspace $(\alpha +\h^\circ)\cap\t^*$, where $\h^\circ$ is the
annihilator of $\h$ in $\fg^*$.  The intersection
$\Phi(M_{prin})\cap\cpwc$ is a connected, relatively open subset of
$S\cap\cpwc$.
\end{enumerate}
\end{lemma}

\begin{proof}\\
\noindent A.  Consider first the case when $\Phi(M_{prin})$ intersects
$\z(\g)^*$, the fixed point set for the coadjoint action of $G$. Let
$m$ be a point in $\Phi \inv (\z(\g)^*)$ and let $\fh$ be its isotropy
Lie algebra. Since the restriction $\Phi|\,M_{prin}$ has constant rank
by Remark~\ref{cons_rank_remark}, there exists a $G$-invariant open
neighborhood $U\subset M_{prin}$ of $m$, such that
$\Phi(U)\subset\g^*$ is a submanifold, and $T_{\Phi(m)}\Phi(U)=
d\Phi_m(T_m\, M)=\h^\circ$.  Since $\Phi(m)$ is fixed by the coadjoint
action, the tangent space $T_{\Phi(m)}\Phi(U)$ is $G$-invariant, which
proves that $\h$ is $G$-invariant. Therefore, the isotropy Lie
algebras of all points in the principal stratum $M_{prin}$ are the
same.  This proves part 2 and 3 of the lemma in the special case, 
and also part 1 for the case $\Phi(m)=0$.
\vskip .1in

\noindent B.  We now reduce the general case to part A, using
symplectic cross-sections.  Given $m\in \Phinv(\cpwc)\cap M_{prin}$,
let $\tau$ be the open wall containing $\Phi(m)$, and $R_\tau$ be the
corresponding natural cross-section.  Denote by $N$ the connected 
component of $R_\tau$ containing $m$. Then $N$ is a Hamiltonian 
$G_\tau$-space, with moment map the restriction $\Phi|N$, 
and $m\in N_{prin}=N\cap M_{prin}$. Since $\Phi(m)\in 
\tau\subset\z(\g_\tau)^*$, part A shows that $\g_m$ is 
an ideal in $\g_{\Phi(m)}$, and that   
every point in $N_{prin}$ has the same isotropy algebra.  
In particular, the isotropy
algebra of points in $\Phinv(\cpwc)\cap M_{prin}$ is locally constant.
We claim that $\Phinv(\cpwc)\cap M_{prin}$ is connected. Indeed,
consider the natural surjective map
$$\pi:\,\Phinv(\cpwc)\cap M_{prin}\ra  M_{prin}/G,\,m\mapsto G\cdot m.$$ 
The fiber of $\pi$ over every point $G\cdot m\in M_{prin}/G$ 
is equal to the intersection 
$G\cdot m\cap \Phinv(\cpwc)\cong G_{\Phi(m)}\cdot m$, 
which is connected since $G_{\Phi(m)}$ is connected. Since 
$M_{prin}$ is connected, it follows that $M_{prin}/G$ is connected. 
Thus $\pi$ has connected target space and connected fibers, and it
follows that $\Phinv(\cpwc)\cap M_{prin}$ is connected. 
This rest of the 
Lemma  is an easy consequence
of the fact that the image of $d\Phi_m$ is the annihilator of $\g_m$.
\end{proof}

\noindent{\it Proof of Theorem \ref{princ cross-section}.}
Let $S\subset \t^*$ be the affine subspace described in Lemma
\ref{lemma stabilizer}. Let $\sigma\subset\cpwc$ be the lowest 
dimensional wall such that $S\cap\cpwc=S\cap\ol{\sigma}$. 
Since the moment map is 
continuous, the moment set $\Phi(M)\cap\cpwc$ is contained in 
the closure of $\Phi(M_{prin})\cap\cpwc$. 
By  Lemma~\ref{lemma stabilizer}, $\Phi(M_{prin})\cap\cpwc$ is an 
open subset of $S\cap\cpwc$. It follows that $\Phi(M_{prin})\cap\sigma$ 
is non-empty, and that the closure of $\sigma$ 
contains the moment set     
$\Phi (M) \cap \cpwc$. It is the smallest wall with this property: 
for any wall $\tau$ with 
$\Phi (M) \cap \cpwc \subset \bar{\tau}$ we have $\sigma \subset
\tau$. 
Consequently $\Phi (M) \cap \sigma = \Phi
(M) \cap U_\sigma$ where $U_\sigma $ is the natural slice
(Remark~\ref{remark natural slice}).  Therefore, by the cross-section theorem,
$$
	Y:= \Phi\inv (\sigma) = \Phi\inv (U_\sigma)
$$
is a symplectic $G_\sigma$-invariant suborbifold of $M$ and $\Phi_Y:=
\Phi|Y : Y \to \sigma \subset U_\sigma$ is a moment map for the action
of $G_\sigma$. Let $\h$ be the stabilizer algebra of points in
$M_{prin}\cap\Phi^{-1}(\cpwc)$. The linear space spanned by $\sigma$
can be identified with the fixed point set $\z(\g_\sigma)^*$ 
for the action of $G_\sigma$ on $\g_\sigma^*$, that is, with the 
annihilator of the
semi-simple part, $[\g_\sigma,\g_\sigma]^\circ\cap\g_\sigma^*$. 
By construction of
$\sigma$, we have $\h^\circ\cap\g_\sigma^*\subset 
[\g_\sigma,\g_\sigma]^\circ\cap\g_\sigma^*$, i.e.
$[\g_\sigma,\g_\sigma]\subset \h$. This shows that the action on $Y$
of the semi-simple part of $G_\sigma$ is trivial, and that the moment
map $\Phi _Y$ is a moment map for the action of the identity component
of the center of $G_\sigma$. That is, the principal cross-section $Y$
is a Hamiltonian {\em torus } orbifold.

We have proved all the assertions of the theorem, except for the fact
that $G\cdot Y$ is dense in $M$ and that $Y$ is connected. 
The
complement to $G \cdot Y \cap M_{prin} = \Phinv(G \cdot \sigma) \cap
M_{prin}$ in $M_{prin}$ is equal to the union of $\Phinv(G\cdot
\tau)\cap M_{prin}$ over all $\tau$ such that $\tau\subset\ol{\sigma}$
and $\tau\not=\sigma$.  By Lemma~\ref{lemma codimension} below, 
these are all submanifolds of 
codimension at least three. It
follows that removing these sets from $M_{prin}$ leaves it dense and
connected.   Hence $G \cdot Y$ is connected and dense in $M$. Since
the quotient map $\g^*\rightarrow \cpwc$ defined by $x\mapsto G\cdot
x\cap\cpwc$ is continuous, $\Phi(Y)$ is dense in $\Phi(M)\cap\cpwc$.
Since $G \cdot Y =G\times_{G_\sigma}Y$
(Remark~\ref{bundle_over_coadj_orbit}) and since the coadjoint orbit
$G/G_\sigma$ is simply connected, it follows that $Y$ is connected.
\hfill$\Box$\\[4pt]
In the proof we used the following Lemma.

\begin{lemma}\label{lemma codimension}
Let $G$ be a compact, connected Lie group, and $M$ a connected
Hamiltonian $G$-orbifold, with moment map $\Phi:\,M\rightarrow\g^*$. Let
$\tau$ be a wall of $\cpwc$ which is different from the principal wall
${\sigma}$.  Then the intersection $\Phinv(G \cdot \tau) \cap
M_{prin}$ is either empty or is a submanifold of codimension at least
3 in $M_{prin}$.
\end{lemma}

\begin{proof}
Let $R_\tau=\Phi^{-1}(U_\tau)$ be the natural cross-section 
corresponding to $\tau$. Let $\z(\g_\tau)$ be the Lie algebra of the 
center of $G_\tau$.
Let $N$ be a connected component of $R_\tau\cap M_{prin}$ such that
$\Phi(N)\cap\tau\not=\emptyset$.  Let $m\in N$, and $\h=\g_m$.  By
Lemma \ref{lemma stabilizer}, 
$\Phi|N$ is a submersion onto an open subset
of the affine space $\g_\tau^*\cap(\Phi(m)+\h^\circ)$, which by assumption
is not contained in $\z(\g_\tau)^*$. The fixed point set for the 
$G_\tau$-action on $\g_\tau^*\cap(\Phi(m)+\h^\circ)$ is the proper 
affine subspace
$\z(\g_\tau)^*\cap(\Phi(m)+\h^\circ)$. 
Since $G_\tau$ is non-abelian, this subspace has 
codimension at least three.  Notice that $\z(\g_\tau)^*$ can be 
identified with the linear space spanned by $\tau$.
Since $\Phi|N$ has constant rank, it 
follows that   
$\Phinv(\tau)\cap M_{prin}$ has codimension greater or equal than 3 in
$R_\tau$. By Remark~\ref{bundle_over_coadj_orbit} this implies that
$$
\Phinv(G\cdot\tau)\cap M_{prin}=G\times_{G_\tau}
(\Phinv(\tau)\cap M_{prin})\subset G\times_{G_\tau}R_\tau
$$ 
has codimension at least 3 in $M_{prin}$.
\end{proof}

\begin{remark}
In the proof of Theorem~\ref{princ cross-section}, we have shown that
the semi-simple part $[\g_\sigma,\g_\sigma]$ of $\g_\sigma$, where
$\sigma$ is the principal wall, is contained in the principal
stabilizer algebra $\h$ for points in $M_{prin}\cap\Phinv(\cpwc)$.  We
also have $\h\subset \g_\sigma$, by equivariance of the moment map.
Therefore, the commutators $[\h,\h]$ and $[\g_\sigma,\g_\sigma]$ are
equal. It follows that the principal stabilizer algebra uniquely
determines the principal wall.
\end{remark}

\section{Hamiltonian torus actions on non-compact orbifolds}

As we have seen, the principal cross-section need not be compact even
if the original manifold is compact.  Thus we need to generalize
Theorem~\ref{AbConv} to include a class of torus actions on
non-compact orbifolds with not necessarily proper moment
maps. Instead, we require that the moment map $\Phi :M\to \ft^*$ is
proper as a map into a convex open set $\sigma$ of $\ft^*$, i.e., that
$\Phi(M) \subset \sigma$ and that for every compact $K \subset \sigma$
the preimage $\Phi^{-1}(K)$ is compact.  This criterion is motivated
by the following fact: if a Lie group acts on a symplectic orbifold
with a proper moment map, then the induced moment map on the principal
cross section is proper as a map into the principal wall.  We extend
Theorem~\ref{AbConv} to this case by using symplectic cuts to
``compactify'' $M$.

First we make an elementary observation. 

\begin{lemma}\label{containing polytope}
Let $(M,\omega)$ be a compact, connected symplectic orbifold, and let
$\Phi: M \to \t^*$ be a moment map for a Hamiltonian torus action on
$M$. Assume that  $\Phi$ is proper as a map into a convex open set
$\sigma\subset \t^*$. 
Then for any compact set $K\subset \sigma$ there
exists a generic rational polytope $P\subset \sigma$, such that $K$ is
contained in the interior of $P$.
\end{lemma}

\begin{proof}
For any point $x\in \sigma$ there exists a polytope $P_x\subset
\sigma$ with  rational vertices which contains $x$ in
the interior.  The collection $\{{\rm int}(P_x): x\in K\}$ is a cover of $K$.
Since $K$ is compact there exists a finite subcover ${\rm int}(P_{x_1})$,
$\ldots$, $   {\rm int}(P_{x_s})$.  Take the convex hull $P$ of the union
$P_{x_1}\cup \ldots \cup P_{x_s}$. If $P$ is generic, it is the
desired polytope. If it is not, perturb it to be generic.
\end{proof}

\begin{theorem} \label{NCabconv} 
Let $(M,\omega)$ be a connected symplectic orbifold, and let
$\Phi: M \to \t^*$ be a moment map for a Hamiltonian torus action on
$M$. If $\Phi$ is proper as a map into a convex open set
$\sigma\subset \t^*$, then
\begin{enumerate}
\item the image $\Phi(M)$ is convex, 
\item each fiber of $\Phi$ is connected, and
\item if for every compact subset $K$ of $\ft^*$, the list 
of isotropy algebras for the $T$-action on  $\Phi\inv (K)$
is finite, then the image $\Phi(M)$ is the intersection of $\sigma$ 
with a rational locally  polyhedral set.
\end{enumerate}
\end{theorem}

\begin{proof}
1. Consider $m_0, m_1 \in M$. Since $M$ is a connected orbifold, there
exists a path $\gamma:\,[0,1]\rightarrow M$ such that $\gamma(0)=m_0$
and $\gamma(1)= m_1$. Since $\gamma ([0 , 1])$ is compact, by
Lemma~\ref{containing polytope} there exists a generic rational 
polytope $P \subset \sigma$ such that $\Phi(\gamma(t))$ is in the
interior of $P$ for all $t\in[0,1]$.  The points $m_0$ and $m_1$ are
contained in the same connected component $N$ of the cut space $M_P$.
Let $\Phi_N:\,N\rightarrow\t^*$ be the induced moment map.  Since
$\Phi_N(N)$ is a convex polytope by Theorem~\ref{AbConv}, it contains
the line segment joining $\Phi(m_0)$ to $\Phi(m_1)$. Since
$\Phi_N(N)\subset\Phi(M)$, this proves that $\Phi(M)$ is convex.\\

\noindent 
2. A similar argument shows that if $\Phi(m_0)=\Phi(m_1)=x$, then $m_0$
and $m_1$ are contained in the same connected component of
$\Phi^{-1}(x)$, since the fibers of $\Phi_N$ are connected by Theorem
\ref{AbConv} and since for points $x$ in the interior of $P$ the
fibers of $\Phi$ and of $\Phi_P: M_P \rightarrow \t^*$ 
are the same (cf. Proposition~\ref{prop_symplectic_cuts}).  
Therefore the fibers of $\Phi$ are connected.\\

\noindent 
3.  For each $x\in\Phi(M)$, let
$$
C_x=\{x+t(y-x)|\,y\in\Phi(M)\quad \hbox {  and } t\ge 0\}
$$ 
be the cone over $\Phi(M)$ with vertex $x$.  Define $A$ to be the
intersection of all the cones $C_x$ for $x\in \Phi (M)$, $A:=
\cap_{x\in \Phi (M)} C_x$.  Since $\Phi(M)$ is convex and is
relatively closed in $\sigma$, Lemma~\ref{convex sets} below (applied
to the closure $X$ of $\Phi(M)$ and to $S=\sigma$) shows that
$\Phi(M)=A\cap\sigma$.
 
We claim that $A$ is the desired rational locally polyhedral set.
First we show that all cones $C_x$ are rational polyhedral cones.  Let
$x\in\Phi(M)$, and let $P\subset\sigma$ be a generic polytope
containing $x$ in its interior.  Since $\Phi(M)$ is convex, $C_x$ is
also the cone over $\Phi_P(M_P)$ with vertex $x$.  In particular,
$C_x$ is a rational polyhedral cone, and the tangent space of each
facet is of the form $\t_i^\circ$, where $\t_i$ is an isotropy Lie
algebra for the action of the torus $T$ in a neighborhood of
$\Phinv(x)$.

Suppose next that $K\subset \t^*$ is a compact convex subset with 
non-empty interior. Since the number of isotropy algebras 
for the $T$-action on $\Phinv(K)$ is finite, 
it follows that up to translation, the list of cones
$C_x$ for $x\in K\cap \Phi(M)$ is finite.  Moreover, if $C_x$ is a
translation of $C_y$, then $C_x=C_y$ since by convexity of $\Phi(M)$
the line segment $\overline{xy}$ is contained in both
$C_x$ and $C_y$.  It follows that the collection of cones $C_x$, $x\in
K\cap \Phi(M)$ itself is a finite list $C_1,\ldots,C_N$. 
This shows that $A$ is a rational locally polyhedral set.
\end{proof}

\begin{lemma}\label{convex sets}
Let $V$ be a vector space, and $X,S\subset V$ convex subsets 
with $X$ closed. For every $x\in X$, let 
$C_x=\{x+t(y-x)\,|\,y\in X,\,t\ge 0\}$.
Then 
$$ X\cap S = \big(\bigcap_{x\in X\cap S}C_x\big)\cap S.$$
\end{lemma}

\begin{proof}
The inclusion ``$\subset$'' is obvious. To prove the opposite
inclusion, assume that $(X\cap S)$ and $S-(X\cap S)$ are non-empty
(since otherwise there is nothing to prove).  Let $y\in S-(X\cap S)$.
We have to show that $y\not\in \bigcap_{x\in X\cap S}C_x$.  Let $r$ be
a ray with vertex $y$ which intersects $X\cap S$ nontrivially.  Since $X$ is
closed and convex $r\cap X$ is either a closed ray or a closed line
segment.  Let $x$ be the point in $r\cap X$ closest to $y$.  Then $y$
does not lie in $C_x$.  Since $S$ is convex, $x$ is in $X\cap S$.
\end{proof}

In the last section, we will use the following corollary to 
Theorem \ref{NCabconv}. 

\begin{cor} \label{NCabmin}  
Let $(M,\omega)$ be a connected symplectic orbifold, and let
$\Phi: M \to \t^*$ be the moment map for a Hamiltonian torus action on
$M$. If $\Phi$ is proper as a map into a convex open set $\sigma$,
then for every $\xi \in  \ft$, 
every local minimum of the function $\Phi^\xi$ is a global minimum,
where $\Phi^\xi(m) := \left< \Phi(m), \xi \right>$.
\end{cor}

\begin{proof}
Since the moment set $\Phi (M)$ is convex, its intersection with the
affine hyperplanes $\{x \in \ft^*\mid \xi (x) = a\}$ is
connected for all $a\in \bbR$. Since the fibers of $\Phi$ are connected
this implies that 
$$ 
(\Phi^\xi)\inv (a) = \Phi\inv (\{x \in
\ft^*\mid \xi (x) = a\}) 
$$ 
is connected for all $a$.  The result follows.  
\end{proof}

\section{Proof of non-abelian convexity}
\label{secproof}

By now, we have done all the necessary hard work.  With a little point
set topology, we bring together the results in the previous two
sections to prove our main theorem.

\noindent
{\bf Theorem~\ref{MainThm} } {\em Let a compact, connected Lie group
$G$ act on a connected symplectic orbifold $(M,\omega)$ with a proper
moment map $\Phi:M \to \fg^*$.
\begin{enumerate}
\item\label{ThmA1} Let $\cpwc$ be a closed Weyl chamber for $G$.  The
{\bf moment set} $\Phi(M) \cap \cpwc$ is a convex  rational locally
polyhedral set. In particular, if $M$ is compact then the moment set
is a convex rational polytope.
\item\label{ThmB1}
Each fiber of the moment map $\Phi$ is connected.
\end{enumerate} }

\begin{proof}
Let $\sigma$ be the principal wall.  By Theorem~\ref{princ
cross-section}, the principal cross section $Y := \Phinv(\sigma)$ is a
connected symplectic orbifold, and the restriction of $\Phi$ to $Y$ is
a moment map for the action of the maximal torus of $Y$.  Since $\Phi
: M \to \fg^*$ is proper, the restriction $\Phi|Y: Y \to \ft^*$ is
proper as a map into the open convex set $\sigma$.  Therefore, by
Theorem~\ref{NCabconv}, the image $\Phi(Y)$ is a convex set and is the
intersection of $\sigma$ with a locally polyhedral set $P$, that is,
$\Phi(Y)= \sigma \cap P$.  By Theorem \ref{princ cross-section},
$\Phi(M)\cap\cpwc=\ol{\Phi(Y)}$. Since the closure of a convex set is
convex, the moment set $\Phi(M)\cap\cpwc$ is convex. Since the closure
of $\sigma $ is a polyhedral cone and since the intersection of the
interior of the locally polyhedral set $P$ with $\sigma$ is nonempty,
the closure of the intersection $\sigma \cap P$ is the intersection
$\bar{\sigma} \cap P$.  Therefore $\Phi(M)\cap\cpwc=\ol{\sigma \cap P}
= \bar{\sigma} \cap P$ is a locally polyhedral set. Since both $P$ and
$\bar{\sigma}$ are rational, the moment set is rational.

It remains to prove that the fiber $\Phi^{-1}(x)$ is connected for all
$x \in \g^*$.  Since the fibers of $\Phi _Y =\Phi |Y $ are connected,
we know that the fibers of the restriction $\Phi | G\cdot Y$ are
connected. Since $G\cdot Y$ is dense in $M$, we would like to conclude
that {\em all} fibers of $\Phi$ are connected.  However, this does not
immediately follow.\footnote{Consider the map 
$f$ from $S^2\subset {\bbR}^3$ to $ S^1\subset {\bbC}$
given by $f(x_1,x_2,x_3) = e^{i \pi  x_3}$.  The map $f$ is proper and
all fibers are connected, except for the fiber over $z=-1$.} \,\,  

First, observe that $\Phinv (x)$ is connected for all $x \in \fg^*$
$\Leftrightarrow$ $\Phinv (G\cdot x)$ is connected for all
$G\cdot x \in \fg^*/G$:
$(\Rightarrow)$ is true since $\Phinv (G\cdot x) = G\cdot \Phinv
(x )$ and the group $G$ is connected;
$(\Leftarrow)$ is true since $\Phinv (G\cdot x)/G = \Phinv
(x )/G_x$ and $G_x$ is connected.

Now, we'll show that the preimage-s of orbits $\Phinv (G\cdot x)$ are
connected. 
For a point $x \in \Phi (M) \cap \cpwc$, let $B$ be an open ball
(with respect to a Weyl group invariant metric) in $\ft^*$ centered at
$x$.  Then $G\cdot (B\cap \cpwc)$ is an open set containing the
orbit $G\cdot x$.  
It enough to show that
the preimage $\Phinv (G\cdot (B\cap \cpwc))$ is connected
(see Lemma~\ref{top_fact}). 
The intersection $\Phinv (G\cdot (B\cap \cpwc)) \cap G\cdot Y$ is
dense in $\Phinv (G\cdot (B\cap \cpwc))$ and $\Phinv (G\cdot (B\cap
\cpwc)) \cap G\cdot Y = G\cdot \Phinv (B\cap \sigma)$, where $\sigma $
is the principal wall.  Since $B\cap \sigma $ is connected and $\Phinv
(y)$ is connected for any $y \in \sigma$, the set $G\cdot
\Phinv (B\cap \sigma)$ is connected.  Since $G\cdot \Phinv (B\cap
\sigma)$ is dense in $\Phinv (G\cdot (B\cap \cpwc))$, the set $\Phinv
(G\cdot (B\cap \cpwc))$ is connected.
This proves that the fibers of the moment map are connected.
\end{proof}
\noindent
In this proof, we used the following topological fact:

\begin{lemma}\label{top_fact}
Let $X$ be a metric space, $f:X \to \bbR^n$ a proper 
continuous map, $\{U_i\}_{i=1}^\infty$ a  sequence of bounded
open sets in $\bbR^n$ with  $U_{i+1} \subset U_{i}$,  $f^{-1} (U_i)$
connected and $C= \bigcap U_i$ nonempty.  Then $f^{-1} (C)$ is also
connected.
\end{lemma}
\begin{proof} Suppose not. Then there are open sets $V$ and $W$ with 
$V\cap W = \emptyset$ and $f\inv (C)\subset V \cup W$.  Since $f^{-1} (U_i)$
is connected there exists $x_i \in f^{-1} (U_i)$ with $x_i\not \in V \cup
W$.  Since the $U_i$'s are bounded, the sequence $f(x_i)$ has a convergent
subsequence and its limit $y$ lies in $C$.  Since $f$ is proper we may
assume, by passing to subsequences, that $f(x_i) \to y$ and that $x_i
\to x$ for some $x \in f^{-1} (y)\subset f^{-1} (C) \subset V \cup
W$. This contradicts the construction of the sequence $\{x_i\}$.
\end{proof}
\begin{remark}
The properness condition on the moment map may be
relaxed to require only that there is an invariant open set $V$
of $\fg^*$ with $\Phi (M) \subset V$, $\Phi :M \to V$ is
proper and and $V\cap \cpwc$ is convex.
\end{remark}

\section{Local moment cones}

A result of Sjamaar \cite{Sj} says that the moment set near a point
$x\in \Delta$ can be read off from local data near an arbitrary point
in the fiber $\Phinv(x)$.  Later Yael Karshon provided a different
proof of Sjamaar's theorem \cite{YK}. In this section, we extend
Sjamaar's results to the orbifold setting, using symplectic cuts.  Our
proof is short and has the advantage that it works in the orbifold
setting.

\begin{theorem}\label{local_cones_thm}
Let $\Phi : M \to \fg^*$ be a proper moment map for an action of a
compact connected Lie group $G$ on a connected symplectic orbifold
$(M, \omega)$.  For every point $m\in M$ and every $G$-invariant
neighborhood $U$ of $m$ in $M$, there exists a $G$-invariant
neighborhood $\calV$ of $\Phi (m)$ in $\fg^*$ such that
$$
\calV \cap \Phi (U) \cap \cpwc =\calV \cap \Phi (M) \cap \cpwc .
$$
\end{theorem}

We will deduce Theorem~\ref{local_cones_thm} from the following
result, which does not require properness of the moment map.

\begin{theorem}\label{local_cones_thm2}
Let  $\Phi : M \to \fg^*$ be a moment map for an action of a
compact connected Lie group $G$ on a symplectic orbifold
$(M, \omega)$. For every point $m\in \cpwc$, there exists a rational 
polyhedral cone $C_m\subset\t^*$ with vertex at $ \Phi(m)$,  
such that for every sufficiently 
small $G$-invariant neighborhood $U$ of $m$, the image $\Phi(U)\cap\cpwc$   
is a neighborhood of the vertex of $C_m$. The cone $C_m$ is called the 
{\bf local moment cone} for $m$.
\end{theorem}

\noindent We first prove Theorem \ref{local_cones_thm}, using 
Theorem \ref{local_cones_thm2}.\vskip .1in 

\begin{proof}
It suffices to check for arbitrarily small $G$-invariant neighborhoods. 
By Theorem~\ref{local_cones_thm2}, there exists an open neighborhood
$\calV$ of $\Phi(M)$ and a rational polyhedral cone 
$C_m \subset \t^*$ such
that $\Phi(U)\cap \cpwc\cap \calV = C_m \cap \calV$.
Since $\Phi(M) \cap \cpwc$ is a convex rational locally polyhedral 
set, there exists a cone $C'_m$ with vertex $\Phi(m)$ such that
$\Phi(M) \cap \calV \cap \cpwc = C'_m \cap \calV$, after shrinking
$\calV$ sufficiently.
Clearly, $C_m \subset C'_m$, so
if $C'_m \neq C_m$, then there exists a
point $x \in \ft^*$ which lies on a face of $C_m$, but in the interior
of $C'_m$. Since $C_m$ is closed and convex there exists $\xi\in\t$
such that a point in $U$ mapping to $x$ is a minimum for $\Phi^\xi$ in
$U\cap Y$, but is not a global minimum for $\Phi^\xi$ on the principal
cross section $Y$.  This contradicts Corollary~\ref{NCabconv}.
Therefore $C_m = C'_m$ and hence 
$$(\Phi(M) \cap\cpwc)\cap \calV=C'_{m}\cap\calV=C_m\cap\calV=
(\Phi(U) \cap \cpwc)\cap \calV.$$
\end{proof}

To prove Theorem~\ref{local_cones_thm2}, we will use the orbifold
version of the local normal form theorem due to Marle \cite{ma} and,
independently, to Guillemin and Sternberg \cite{g-s:normal}. The
theorem asserts that an invariant neighborhood of a point $m$ in a
Hamiltonian $G$-manifold $M$ is completely determined (up to
equivariant symplectomorphism) by two pieces of data: (1) the value of
the moment map $\Phi$ at $m$ and (2) the symplectic slice
representation of the isotropy group $G_m$ of the point.  Recall that
the {\bf symplectic slice} at a point $m$ is the largest symplectic
subspace in the fiber at $m$ of the normal bundle to the orbit $G\cdot
m$.  The analogous result holds for Hamiltonian $G$-orbifolds (cf.\
\cite{LT}) --- the only difference is that the symplectic slice is no
longer a vector space.  Instead it is the quotient of a vector space
by a linear action of a finite group.

\begin{lemma}\label{locsymp}  
Let $G$ be a compact Lie group, let $(M,\omega)$ be a Hamiltonian
$G$-orbifold, with moment map $\Phi$, and $m\in\Phinv(0)$ a point in
the zero level set. Let $\Gamma$ be the orbifold structure group of
$m$ and $G_m\times V/\Gamma\ra V/\Gamma$ be the symplectic slice
representation at $m$. For every $G_m$-equivariant splitting
$\fg^*\cong\fg_m^\circ \oplus\fg^*_m$, there is a $G$-invariant
symplectic form on the orbifold $F=G\times_{G_m} \left(\fg_m^\circ
\times V/\Gamma \right) $ such that the moment map $\Phi_F : F \to
\fg^*$ for the $G$-action is given by
$$ 
\Phi _F ([g, \eta, v])= g\cdot \left(\eta + \phi_{V/\Gamma}
(v)\right),
$$ 
where $\phi_{V/\Gamma} : V/\Gamma \to \fg_m^*$ is the moment map for the
slice representation of $G_m$.  Moreover, the embedding of the orbit
$G\cdot m$ into $F$ as the zero section is isotropic and the
symplectic slice at $[1,0,0]$ is $V/\Gamma$.  Consequently there
exists a $G$-equivariant symplectomorphism $\lambda$ from a
neighborhood $U$ of $G\cdot m$ in $M$ to a neighborhood $U'$ of the
zero section in $F$, and $\lambda^*\Phi_F=\Phi$ over $U$.
\end{lemma}

We are now ready to prove Theorem~\ref{local_cones_thm2}. 

\begin{proof}
1. We begin by considering the special case $\Phi(m) = 0$.
Because we are considering arbitrarily small neighborhoods,
by Lemma~\ref{locsymp}  
it  suffices to consider neighborhoods of the zero section
of the model space $F$.
Let $V/\Gamma$ be the symplectic slice at $m$. Choose a $G_m$-invariant 
complex structure on $V/\Gamma$ which is compatible with the 
symplectic form; let $\rho$
denote the norm squared of the induced metric.  The model $F$ is a
complex orbi-bundle over $G\times _{G_m}\fg_m^\circ$. The
multiplication action of $U(1)$ on the fibers of $F \to G\times
_{G_m}\fg_m^\circ$ is a Hamiltonian action on the orbifold $F$, with
moment map being the function $\rho$ defined above. The action of
$U(1)$ commutes with the action of $G$.  Consequently the moment map
$\hat{\Phi}_F : F \to \bbR \times \g^*$ for the action of $U(1)
\times G$ on $F$ is given by
$$ 
\hat{\Phi}_F ([g, \eta, v])= 
\Big( \rho([g, \eta, v]),  g\cdot \big(\eta + \phi_{V/\Gamma} (v)\big) \Big).
$$
Since $\hat{\Phi}_F$ is proper, we can apply Theorem~\ref{MainThm}.
Therefore, $\hat{\Phi}_F(F)\cap (\bbR \times \cpwc)$ is a convex
rational locally polyhedral set. In fact, since $\hat{\Phi}_F$ is homogeneous
(i.e., equivariant with respect to the action of $\bbR_+$ which on $F$
is given by $t\cdot[g,\eta, v] = [g, t\eta, \sqrt {t} v]$ and on
$\fg^*$ by  multiplication) and since the number of orbit types in $F$
is finite, $\hat{\Phi}_F(F)\cap (\bbR \times \cpwc)$ is a convex
rational polyhedral cone.  Since $\Phi _F (F)\cap \cpwc$ is the image of
$\hat{\Phi}_F(F)\cap (\bbR \times \cpwc)$ under $\pi:\bbR \times \cpwc
\to \cpwc$, it is also a convex rational polyhedral cone. \vskip .1in

2. Choose a $G$-invariant metric on $\g^*$, and let $\hat{\rho}:\,F\ra
\bbR$ be defined by $\hat{\rho}([g,\eta,v])= \rho(v)+||\eta||$. Then
$\hat{\rho}$ is homogeneous and $\hat{\rho}^{-1}(0)$ is the zero
section of $F$.  Choose $\epsilon>0$ sufficiently small so that the
set $\{x \in F\mid \hat{\rho}(x) <\epsilon\}$ is contained in $U$.
Let $\lambda$ be the lower bound on $\hat{\rho}^{-1}(\epsilon)$ of the
function
$$
f:F \to \bbR,\, f(x)=\rho(x)+||\Phi_F(x)||
$$ 
on $F$.  Since $f\ge\hat{\rho}$, we have $\lambda\ge\epsilon$. Let
$\hat{\calV}:=\{(t,\alpha)\in \bbR\times\g^*|\,\,|t|+||\alpha||
<\lambda\}$.  Since $||\hat{\Phi}_F(x)||\ge \lambda$ if
$\hat{\rho}(x)\ge \epsilon$, one has $\hat{\Phi}_F(U)\supset
\hat{\Phi}_F(F)\cap \hat{\calV}$, which is a $G$-invariant
neighborhood of $0\in \hat{\Phi}_F(F)$.  Since $\hat{\Phi}_F(F)\cap
(\bbR \times \cpwc)$ is a convex rational polyhedral cone, the image
of every neighborhood of the vertex of this cone under the projection
$\pi:\bbR \times \t^* \to \t^*$ is a neighborhood of the vertex of
$\Phi_F(F)\cap\cpwc= \pi(\hat{\Phi}_F(F)\cap \bbR\times\cpwc)$. It
follows that $G\cdot \pi(\hat{\calV}\cap (\bbR\times \cpwc))$ is a
$G$-invariant neighborhood of $0$ in $\fg^*$ with the required
property.
\vskip .1in

3. Finally, we use cross sections to reduce the general case to  
the case that $\Phi(m)=0$.  Let $\tau$ be the (open) wall
containing $\Phi (m)$. Let $U_\tau$ be the corresponding natural slice
and $R_\tau = \Phi\inv (U_\tau)$ the natural cross-section.
By equivariance of the moment map, we have
$$
  \Phi (G\cdot R_\tau) \cap \cpwc =  \Phi (R_\tau )\cap \cpwc .
$$
The Cartan subalgebras of $G$ and $G_\tau$ are equal, but a Weyl
chamber $(\ft^*_\tau)_+ $ of $G_\tau$ is the union of certain Weyl
chambers of $G$.  Nonetheless, $\Phi (R_\tau) \cap (\ft^*_\tau)_+ =
\Phi (R_\tau) \cap \cpwc$, which means that we can apply Remark
\ref{general argument} to reduce to the case $\Phi (m) =0$.
\end{proof}

\begin{remark}
A corollary to the above Theorems (or rather their proofs) 
is that if $M$ is a Hamiltonian $G$-orbifold, with moment map 
$\Phi$, the local moment cone $C_m$ for a point $m\in\Phinv(\cpwc)$ 
does not vary as $m$ varies in a connected component of a fixed fiber 
of $\Phi$.
Indeed, if $M$ is connected and $\Phi$ is proper this follows 
from Theorem \ref{local_cones_thm} and Theorem 
\ref{local_cones_thm2}, because the local moment 
cone is equal to the cone with vertex at $\Phi(m)$ over the 
moment set $\Delta$. 
For the general case, we note that as above, we can use cross-sections
to reduce to the case $\Phi(m)=0$, and since the statement is 
local, it is enough to check it for the local normal form, $F$. 
Now the moment map $\Phi_F$ is not proper. However, we may replace 
$F$ by its cut $F_{(-\infty,\epsilon]}$ with respect to the 
$S^1$-action generated by $\rho$. Since the moment map for the 
cut space is proper, the claim follows.
\end{remark}


\begin{thebibliography}{MMMM}
\bibitem[A]{A} M. F. Atiyah, Convexity and commuting hamiltonians, {\it
Bull. London Math. Soc.} {\bf 14} (1982), 1-15.

\bibitem[B]{B} M. Brion, { Sur l'image de l'application moment}, in
{\em S\'eminaire d'alg\`ebre Paul Dubreuil et Marie-Paule Malliavin}
(Paris, 1986) (M.-P. Malliavin, ed.), Lecture Notes in Mathematics,
vol 1296, Springer-Verlag, Berlin-Heidelberg-New York, 1987, 177--192.

\bibitem[CDM]{CDM} M. Condevaux, P. Dazord, P. Molino, 
G\'eom\'etry du Moment, in {\em Travaux du S\'eminar Sud-Rhodanien I},
Publ.\ D\'ep.\ Math., Univ.\ Claude Bernard - Lyon 1 (1988), 131--160.

\bibitem[GLS]{GLS} V. Guillemin, E. Lerman and
 S. Sternberg, {\em Symplectic fibrations and multiplicity diagrams},
 Cambridge University Press (to appear).

\bibitem[GS]{g-s:normal}
V. Guillemin and S. Sternberg,
 A normal form for the moment map, in: {\em Differential Geometric Methods
in Mathematical Physics} (S. Sternberg, ed.),  Reidel Publishing Company,
Dordrecht, 1984.

\bibitem[GS1]{GSc} V. Guillemin and S. Sternberg, Convexity properties of
the moment mapping I and II, {\it Invent.\ math.} {\bf 67} (1982),
491--513,  {\it Invent. math.} {\bf 77}  (1984), 533--546.

\bibitem[GS2]{GS} V. Guillemin and S. Sternberg, {\em Symplectic techniques 
in physics}, Cambridge University press, 1990.

\bibitem[H]{Hec} G. Heckman, Projection of orbits and asymptotic
behavior of multiplicities for compact Lie groups, Ph.D. thesis,
Leiden, 1980.

\bibitem[HNP]{Hilgert} J. Hilgert, K.-H. Neeb and W. Plank, Symplectic
convexity theorems and coadjoint orbits, {\em Comp.\ Math.} {\bf
94} (1994), 129--180.

\bibitem[Ka]{YK} Y. Karshon, unpublished notes.

\bibitem[Ki1]{Kbook} F. Kirwan, {\it Cohomology of Quotients in Symplectic
and Algebraic Geometry}, Princeton University Press, Princeton, 1984.

\bibitem[Ki2]{Kconvex} F. Kirwan, Convexity properties of the moment
mapping III, {\it Invent.\ math.} {\bf 77} (1984), 547--552.

\bibitem[Ko]{Kostant}  B. Kostant, On convexity, the Weyl group, and the
Iwasawa decomposition, {\it Ann. Sci. Ec. Norm. Sup.} {\bf 6} (1973), 413--455.

\bibitem[L]{L} E. Lerman, Symplectic cuts, {\em Math.\ Research Lett.}
{\bf 2} (1995),  247--258.

\bibitem[LT]{LT} E. Lerman and S. Tolman, Symplectic toric orbifolds,
{\tt dg-ga/9412005}; Hamiltonian torus actions on
symplectic orbifolds and toric varieties {\tt dg-ga//9511008}.


\bibitem[Ma]{ma}
C. M. Marle, { Mod\`{e}le d'action hamiltonienne d'un groupe de
Lie sur une vari\'{e}t\'{e} symplectique}, {\em Rendiconti del
Seminario Matematico} {\bf 43} (1985), 227--251, Universit\`a e
Politechnico, Torino.

\bibitem[M]{Mein}  E. Meinrenken,  Symplectic surgery and the Spin$_c$-Dirac
operator, to appear in {\em Adv.\ Math.},
{\tt  dg-ga/9504002}.  

\bibitem[NM]{NM} L. Ness, A stratification of the null cone via the
moment map, {\em Am. J. of Math.} {\bf 106}, (1984), 1281-1329, with an
appendix by D. Mumford.
\bibitem[Sa]{Sa}
I. Satake, The Gauss-Bonnet theorem for V-manifolds {\em J. Math.\
Soc.\ Japan} {\bf 9} (1957), 464--492. 


\bibitem[Sj]{Sj} R. Sjamaar, Convexity properties of the moment map
re-examined, preprint, to appear in {\em Adv.\ Math.}, 
{\tt dg-ga/9408001}.  

\bibitem[W]{W} C. Woodward, The classification of transversal
multiplicity-free group actions, {\em Ann.\ Global Anal.\ Geom.} {\bf
14} (1996), 3--42.

\end{thebibliography}
\end{document}